\documentclass[12pt]{article}
\usepackage{epsfig}
\usepackage{epsf}
\usepackage{pifont}

\newcommand{\R}{{{\cal R}}}
\newcommand{\N}{{{\cal N}}}

\newcommand{\LL}{{\cal L}}
\newcommand{\FF}{{{\cal F}}}
\newcommand{\G}{{{\cal G}}}
\newcommand{\Ob}{{{\cal O}}}
\newcommand{\Ez}{{{\rm E}_0}}
\newtheorem{Theorem}{Theorem}

\newtheorem{Lemma}{Lemma}
\newtheorem{Definition}{Definition}
\begin{document}

\title{Likelihood For Generally Coarsened Observations From Multi-State Or Counting Process Models}
\author{ Daniel Commenges, Anne G\'egout-Petit\\ INSERM; Universit\'e Victor Segalen Bordeaux 2}
\date{}
\maketitle

{\bf Running title: Likelihood for general coarsening}

\bigskip

    {\bf ABSTRACT. We consider first the mixed discrete-continuous scheme of observation  in  multi-state models; this is a classical pattern in epidemiology because very often clinical status is assessed at discrete visit times while times of death or other events are observed exactly. A heuristic likelihood can be written for such models, at least for Markov models; however a formal proof is not easy and has not been given yet. We present a general class of possibly non-Markov multi-state models which can be represented naturally as  multivariate counting processes. We give a rigorous derivation of the likelihood  based on applying Jacod's formula for the full likelihood and taking conditional expectation for the observed likelihood. A local description of the likelihood allows us to extend the result to a more general coarsening observation scheme proposed by Commenges \& G\'egout-Petit (2005). The approach is illustrated by considering  models for dementia, institutionalization and death.}
\vspace{2mm}


\section{Introduction}
Multi-state models have been proposed for a long time, in particular in biological applications. Until the seventies however the attention was essentially focused on homogeneous Markov models. Non-homogeneous Markov models were studied by Fleming (1978), Aalen \& Johansen (1978); semi-Markov models were also considered by Lagakos, Sommer \& Zelen (1978); the most studied model was the illness-death model (Andersen, 1988; Keiding, 1991). A thorough account of these developments can be found in Andersen {\em et al.} (1993) where the counting process theory is used to obtain rigorous results in this field; see also Hougaard (2000).

Another stream of research was started by Peto (1973) \& Turnbull (1976) who tackled the problem of interval-censored observations in survival data analysis and gave the non-parametric maximum likelihood estimator;  Frydman (1995a, 1995b) extended this issue to the illness-death model in which transition toward illness could be interval-censored while time of death was exactly observed, while Wong (1999) studied a general case of multivariate interval-censored data. The penalized likelihood approach to this problem was proposed by Joly \& Commenges (1999) with an application to AIDS, and by Joly {\em et al.} (2002) and Commenges {\em et al.} (2004) with application to Alzheimer's disease.

In the continuous time
observation scheme it has been shown (Aalen, 1978; Borgan, 1984; Andersen {\em et al.},
1993) that the likelihood could be derived
from Jacod's formula (Jacod, 1975); also in this case the martingale theory yields simple and natural estimators (Andersen {\em et al.}, 1993; Aalen {\em et al.}, 2004). In the case where some transitions
are observed in discrete time while others are observed in
continuous time, that we call the Mixed Discrete-Continuous
Observation (MDCO) scheme, the natural martingale estimators are no longer available, so one has to return to likelihood-based methods. The above cited papers considering this case used
heuristic likelihoods. Commenges (2003) derived the likelihood for
this observation scheme in a Markov illness-death model. The main
aim of this paper is to rigorously derive the likelihood in a more general
framework: i) we consider a general class of multi-state models which may have
any number of states and are not necessarily Markov; ii) we first consider
 the MDCO scheme and we extend to the so-called GCMP (General Coarsening Model for Processes) proposed by Commenges \& G\'egout-Petit (2005). In most of this work we assume that the mechanism leading to incomplete data is ignorable (Gill {\em et al.}, 1997); however we prove ignorability in the special case when death produces a stochastic censoring of the other processes of the model.
Our approach starts by remarking
that most of the useful multi-state models can be directly
formulated in terms of multivariate counting process in a very
natural way, an idea close to the ``composable processes'' studied by Schweder (1970). Then Jacod's formula can be applied to find the
likelihood for continuous time observation. In the case where one
or several of these processes is observed in discrete time the
likelihood can be computed by taking a conditional expectation of
the full likelihood.

In section 2 the heuristic likelihoods for diverse observation
schemes are recalled. Section 3 develops a natural correspondence
between the class of irreversible multi-state processes and multivariate one-jump counting processes. In
section 4 the multivariate counting process representation is
exploited to use Jacod's formula for finding the likelihood in the
mixed discrete-continuous observation scheme; then the result is extended to the GCMP scheme. This general
modeling approach is illustrated in section 5 for describing a
joint model for dementia, institutionalization and death,
presented as a five-state model in Commenges \& Joly (2004) and
showing the benefits of the proposed approach to this case. Section 6 briefly concludes.

\section{Heuristic likelihoods for multi-state models}
\subsection{Notation} A multi-state process $X=(X_t)$ is a
right-continuous process which can take a finite number of values
$\{0, 1, \ldots, K-1\}$. The theory of Markov multi-state models (or
Markov chain models) is well established. The law of a Markov
multi-state process is defined by  the transition probabilities
between states $h$ and $j$ that  we will denote by
$p_{hj}(s,t)=P(X_t=j|X_s=h)$; transition intensities
$\alpha_{hj}(t)$, for $h\ne j$, may be defined when the $p_{hj}(s,t)$'s
are continuous both in $s$ and in $t$ for all $s$ and $t$, as the
following limits (if they  exist):
\begin{equation}
\label{transitions}{\rm lim}_{\Delta t \rightarrow 0}
p_{hj}(t,t+\Delta t)/\Delta t. \end{equation} It  is reasonable in
most applications in epidemiology to think that these limits
exist; it is even reasonable to expect continuous and smooth
transition intensities. We define  $\alpha_{hh}(t)=-\sum _{j\ne h}
\alpha_{hj}(t)$ and $\sum _{j\ne h} \alpha_{hj}(t)$ is the hazard
function associated with the distribution of the sojourn time in
state $h$. When the $\alpha_{hj}(t)$'s do not depend on $t$, the
Markov chain is said to be homogeneous. Transition probabilities
and transition intensities are linked by the forward Kolmogorov
differential equations, so that solving these equations transition
probabilities can be expressed as a function of transition
intensities; this solution can take the form of the product
integral (see Andersen {\em et al.}, 1993). For non-Markov multi-state
models one could define analogously transition probabilities
$p_{hj}(s,t)=P(X_t=j|X_s=h,{\cal F}_{s-})$, where ${\cal F}_{s-}$
is the history before $s$; similarly, transition intensities
$\alpha_{hj}(t;{\cal F}_{t-})$ could be defined. However we are on
a less firm ground because, contrarily to the Markov case, these
quantities are random and the Kolmogorov equations have been given
only for Markov processes.

In the remaining of this paper we shall consider observations for different schemes of observation, leading to incomplete data.
We assume that the mechanism leading to incomplete data is ignorable: this means that we make a correct inference by using the likelihood as if these schemes were deterministic. In particular we can consider that the different times like $C$ and $v_l, l=0,\ldots, m$  involved in this mechanisms, and which are defined below, are fixed. This raises however a problem which will be discussed and solved in section 4.5.

\subsection{Likelihood for continuous
time observations in Markov models} Consider the case where
process $X$ is continuously observed from $v_0$ to $C$.  We
observe that transitions have occurred at (exactly) times
$T_{(1)}< T_{(2)}<\ldots <T_{(M)}$. With the convention
$T_{(0)}=v_0$, the value of the likelihood, conditional on
$X_{v_0}$, on the event $\{M=m\}$ and $\{X_{T_{(r)}}=x_r,
r=0,\ldots,m\}$ is: $${\cal L}= [\prod _{r=1}^{m}
p_{\scriptscriptstyle x_{r-1},x_{r-1}}{\scriptstyle
(T_{(r-1)},T_{(r)}-)}\alpha_{\scriptscriptstyle
x_{r-1},x_{r}}{\scriptstyle (T_{(r)})}]p_{\scriptscriptstyle
x_m,x_m}{\scriptstyle (T_{(m)},C)}.$$ The probability that no
transition happens between $T_{(r-1)}$ and $T_{(r)}-$ given
$X_{T_{(r-1)}}=x_{r-1}$ can easily be computed as
$$p_{\scriptscriptstyle x_{r-1},x_{r-1}}{\scriptstyle
(T_{(r-1)},T_{(r)}-)}= \exp \int_{T_{(r-1)}} ^{T_{(r)}}
\alpha_{\scriptscriptstyle x_{r-1},x_{r-1}}(u) du.$$

\subsection{Likelihood for discrete
time observations in Markov models}
 Consider now the case where $X$ is observed at discrete
times $v_0, v_1,\ldots,v_m$. In this case, we observe a vector of
random variables $(X_{v_0},\ldots,X_{v_m})$ and it is easy to
derive the value of the likelihood, conditional on $X_{v_0}$, on
the event $\{X_{v_{r}}=x_r, r=0,\ldots,m\}$:
$${\cal L}= \prod
_{r=0}^{m-1} p_{x_r,x_{r+1}}(v_r,v_{r+1}).$$

\subsection{Likelihood for mixed discrete-continuous time
observations in Markov models} The most common case in
applications is that some transitions are observed in discrete
times and others in continuous time. A classical example is the
irreversible illness-death model, a model with the three states ``health'',
``illness'', ``death'' respectively labeled $0, 1, 2$; it is
often the case that transition toward the illness state is
observed in discrete time while transitions toward death is
observed in continuous time (Frydman, 1995a; Joly {\em et al.}, 2002).
Let us call $\tilde T$ the follow-up time that is $\tilde T= {\rm
min}(T,C)$, where $T$ is the time of death; we observe $\tilde T$
and $\delta=1_{\{T\le C\}}$. If the subject starts in state
``health'', has never been observed in the ``illness'' state and
was last seen at visit $M$ (at time $v_M$) the likelihood,
conditional on $X_{v_0}=0$, on the event $\{M=m\}$ is:
\begin{equation} \label{likh} {\cal
L}=p_{00}(v_0,v_m)[p_{00}(v_m,\tilde T)\alpha_{02}(\tilde
T)^{\delta}+p_{01}(v_m,\tilde T)\alpha_{12}(\tilde
T)^{\delta}];\end{equation}
if the subject has been observed in
the illness state for the first time at $v_l$ then the likelihood
is:
\begin{equation} \label{liki} {\cal
L}=p_{00}(v_0,v_{l-1})p_{01}(v_{l-1},v_l)p_{11}(v_l,\tilde
T)\alpha_{12}(\tilde T)^{\delta}.\end{equation} These formulae can
be extended rather easily to models with more than three states in
the case where there is one absorbing state, and transitions
toward the absorbing state are observed in continuous time while
transitions toward other states are observed in discrete time
(Commenges, 2002). In the case where transitions toward one state
are observed in discrete time while transitions toward other
states are observed in continuous time the formula becomes
cumbersome: see the example of the Dementia-Institution-Death
model of Commenges \& Joly (2004). These formulae are heuristic.
The same type of formulae can be written for non-Markov model
using random versions of the transition probabilities and
intensities and hoping to be able to compute transition
probabilities in term of transition intensities; see Joly \&
Commenges (1999) for an example of a semi-Markov model.

\section{Multi-State models as multivariate counting processes}
The main motivation for representing multi-state models as
multivariate counting processes is the availability of a formula
giving the likelihood ratio for such processes (not necessarily
Markov) observed in continuous time on $[0,C]$ (Jacod, 1975). It
will be shown in the next section that this can be the basis of a
rigorous derivation of the likelihood also in the MDCO and in the
more general GCMP schemes. It was shown (Borgan, 1984; Andersen {\em et
al.}, 1993) that a multi-state model $(X, \mathbf{\alpha})$, where
$\alpha=(\alpha_{hj}(.); 0 \le h,j \le K-1)$, can be represented
by the multivariate counting process $(N, \mathbf{\lambda})$ where
$N=(N_{hj}, h\neq j, 0\leq h,j \leq K-1)$ and the $N_{hj}$'s count
the number of transitions from state $h$ to state $j$; the
intensity of $N$ is $\mathbf{\lambda}=(\lambda_{hj}, h\neq j,  0
\le h,j \le K-1)$ and $\lambda_{hj}=Y_h(t)\alpha_{hj}(t)$, where
$Y_h(t)=1_{\{X_{t-}=h\}}$.

For irreversible multi-state models a more parsimonious representation is possible.
It is often possible to formulate an epidemiological problem
directly in term of counting processes rather than using the
counting process representation as a mathematical device. Most
multi-state models are in fact used for jointly modeling several
events. For instance the illness-death model is used to jointly
model onset of disease and death. So we can directly model the
problem by considering a bivariate counting process $N=(N_1,N_2)$,
where $N_1$ counts the onset of the disease and $N_2$ counts the
occurrence of death. It can be seen that the multi-state process
can be retrieved by $X_t=\min(2N_{2t}+N_{1t},2)$. The processes
$N$ and $X$ generate the same filtration. Note however that for a
fixed $t$, the random variables $X_t$ and $N_t$ do not generate
the same $\sigma$-field in general because the event
$\{X_t=2\}$ (subject in the ``death'' state at $t$) is the same as $\{N_{2t}=1\}$; however if we know $N_t$, we also know $N_{1t}$, that is we know whether the subject passed through the state ``illness''  before $t$ or not. Note
that the representation based on the basic events of interest is
more economical than the Borgan representation (a bivariate rather
than three-variate process). Similarly the
Dementia-Institution-Death model, which is a five-state model (see
section 5), can be represented by a three-variate process counting
onset of dementia, institutionalization and death; in this model
there are eight possible transitions so that the Borgan
representation would entail a eight-variate process.

 Let us consider the more general problem of jointly modeling the onset of
$p$ types of events, each type occurring just once. This can be
represented by a p-variate counting process $N$, each $N_j$ making
at most one jump. It is possible to construct  a multi-state
process $W$ such that $W_t$ generates the same $\sigma$-field as
$N_t$. A possibility is:
$W_t=N_{pt}2^{p-1}+N_{p-1,t}2^{p-2}+\ldots N_{1t}$, that we denote by
$W_t=N_{pt}N_{p-1,t}\ldots N_{1t}$ (this is the representation of
$W_t$ in base 2); $W_t$ can take $2^p$ integer values in the set
$\{0,1,\ldots,2^p-1\}$. Consider the important case where $N_p$
counts death; it is common in the multi-state representation to
consider that deceased subjects are in the same state, that is we
may construct a more compact multi-state model $X$ defined by
$X_t=\min (W_t, 2^{p-1})$; $X_t$ can take $2^{p-1}+1$ values. For
example the Dementia-Institution-Death model (where $p=3$) has
five states rather than eight. We have exactly the same number of
non-zero transition intensities for $W$ and $X$ and they are
equal; we simply have to rename them. More specifically we have
$\alpha^W_{hj}(.)= \alpha_{hj}(.)$ for $0 \le h<j\le 2^{p-1}$ and
$\alpha^W_{hj}(.)=\alpha_{h2^{p-1}}(.)$ for $0 \le h<2^{p-1}\le j$.

\begin{Theorem} Let $N=(N_1,\ldots,N_p)$ be a counting process
with $N_{jt}\le 1; j=1,\ldots, p;  t\ge 0$ and $p>1$. Consider the
multi-state process $W=(W_t)$ defined by $W_t=N_{pt}\ldots N_{1t}$
in base 2. If, in a given probability measure, $W$ is Markov with
continuous transition intensities $\alpha^W_{hj}(.); 0\le h,j \le
2^p-1$, the generating counting process $N_j, j=1\ldots p$, have
intensities given by:
\begin{equation}
\label{intens-intens}
\lambda_j(t)=1_{\{N_{jt-}=0\}}\sum_{k_1=0}^{1}\ldots
\sum_{k_{p}=0}^{1}
\prod_{l=1}^{p}1_{\{N_{lt-}=k_l\}}\alpha^W_{k_p\ldots
k_{j+1}0k_{j-1} \ldots k_1,k_p \ldots k_{j+1} 1 k_{j-1}\ldots
k_1}(t); t\ge 0 ,\end{equation}
where $k_p\ldots k_{j+1}0k_{j-1}
\ldots k_1$ and $k_p\ldots k_{j+1}1k_{j-1} \ldots k_1$ are base 2
representations of integers. \end{Theorem}

\noindent {\em Proof.}  Lemma 3.3 of Aalen (1978) gives an expression of the c\`{a}dl\`{a}g modification
$(\lambda_j(t+))$ of the c\`{a}gl\`{a}d process $ (\lambda_j(t))$ by :
\begin{eqnarray*} \hspace{-20mm} \lambda_j(t+) & = & \lim _{\delta
\downarrow 0} \frac {1}{\delta}  P[N_{j(t+\delta)}-N_{jt}=1|\FF_t]
\\ & = & \lim _{\delta \downarrow 0}\frac {1}{\delta} P[\cap_{l\ne
j} \{N_{l(t+\delta)}=N_{lt}\}\cap
\{N_{j(t+\delta)}-N_{jt}=1\}|\FF_t]\\ & = & 1_{\{N_{jt}=0\}}\lim
_{\delta \downarrow 0}\frac {1}{\delta} P[\cap_{l\ne j}
\{N_{l(t+\delta)}=N_{lt}\}\cap \{N_{j(t+\delta)}=1\}|\FF_t]\\ & =
& 1_{\{N_{jt}=0\}}\lim _{\delta \downarrow
0}\sum_{k_1=0}^{1}\ldots \sum_{k_{p}=0}^{1}
\prod_{l=1}^{p}1_{\{N_{lt}=k_l\}}\frac {1}{\delta} P[\cap_{l\ne j}
\{N_{l(t+\delta)}=N_{lt}\}\cap \{N_{j(t+\delta)}=1\}|\FF_t]\\ & =
& 1_{\{N_{jt}=0\}}\sum_{k_1=0}^{1}\ldots \sum_{k_{p}=0}^{1}
\prod_{l=1}^{p}1_{\{N_{lt}=k_l\}}\lim _{\delta \downarrow 0}\frac
{1}{\delta} P[W_{(t+\delta)}=k_p \ldots 1 \ldots k_1|W_{t}=k_p
\ldots  0\ldots k_1]\\ & = &
1_{\{N_{jt}=0\}}\sum_{k_1=0}^{1}\ldots \sum_{k_{p}=0}^{1}
\prod_{l=1}^{p}1_{\{N_{lt}=k_l\}}\alpha^W_{k_p\ldots
k_{j+1}0k_{j-1} \ldots k_1,k_p \ldots k_{j+1} 1 k_{j-1}\ldots
k_1}(t), \end{eqnarray*} from which the theorem follows.

A simple way of reading formula (\ref {intens-intens}) is to say
that on an event such that a jump of $N_j$ implies a jump of $W$
from $h$ to $h'$, then $\lambda_j(t)=\alpha^W_{hh'}(t)$. It is
easy to apply this theorem for finding the intensities of the
$N_j$'s as a function of the transition intensities of $X$ defined
as above by first finding the intensities in term of the transition intensities of $W$
and then renaming the transition intensities.

{\em Example: The illness-death model.} Consider a bivariate counting process
$N=(N_1,N_2)$, where $N_1$ counts the onset of the disease and
$N_2$ counts the occurrence of death; thus we have $p=2$.  We can
form the process $(W_t)$ defined by $W_t=N_{2t}N_{1t}$. This
process has four states but states $2$ and $3$ have the same
biological meaning (``dead'') so that they are generally grouped
to obtain the conventional illness-death model $(X_t)$ defined by
$X_t=\min (W_t,2)$. If $W$ is Markov, applying formula
(\ref{intens-intens}) yields  for the intensities of $N_1$ and
$N_2$ respectively: $$\lambda
_1(t)=1_{\{N_{1t-}=0\}}1_{\{N_{2t-}=0\}}\alpha^W_{01}(t)$$ $$
\lambda
_2(t)=1_{\{N_{2t-}=0\}}[1_{\{N_{1t-}=0\}}\alpha^W_{02}(t)+1_{\{N_{1t-}=1\}}\alpha^W_{13}(t)].
$$ In term of the transition intensities of $X$ we obtain:
\begin{equation}\label{int1} \lambda
_1(t)=1_{\{N_{1t-}=0\}}1_{\{N_{2t-}=0\}}\alpha_{01}(t)\end{equation}
\begin{equation} \label{int2} \lambda
_2(t)=1_{\{N_{2t-}=0\}}[1_{\{N_{1t-}=0\}}\alpha_{02}(t)+1_{\{N_{1t-}=1\}}\alpha_{12}(t)].
\end{equation}

 The same approach can be applied if some counting
processes take more than two values even if the trick of the
representation of the state value in base 2 does not work. A
process with more than one jump can represent progression of a
disease: for instance the progressive model of Joly \& Commenges
(1999) can be represented by a counting process making a jump at
HIV infection and another jump at onset of AIDS. Alternatively we can represent progression of the disease
 by two $0-1$ counting processes, one counting HIV infection, the other counting onset of AIDS, with the intensity of the latter being equal to zero if the former has not jumped.  The 10-state
model proposed by Alioum {\em et al.} (2005) can be represented by a
process representing progression of the HIV disease which can make
two jumps and three $0-1$-processes counting HIV diagnosis, inclusion in
a cohort and death or alternatively by five $0-1$-processes if we represent progression of the disease by the two $0-1$ counting processes as discussed above. A problem appears if we make the base-2 construction for such $0-1$ processes. The construction leads to ``phantom'' states which are not relevant, for instance the state AIDS without HIV infection; however we can still represent the relevant multi-state model by putting the transition leading to this state uniformly equal to zero. Theorem 1 can then still be applied.

We can formalize the relationship between the class of irreversible multi-state (IM) processes and one-jump counting (OJC) processes by saying that they are equivalent. This concept of equivalence between classes of processes is based on the following  equivalence relation between processes.

\begin{Definition} Two processes are informationally equivalent if they generate the same filtration.
\end{Definition}

 For each IM  process we can find at least one informationally equivalent OJC process; this is obvious from the Borgan representation. Inversely for each OJC process we can find an informationally equivalent IM process (using the base 2 representation). We can also define canonical processes, that is simple representants of an equivalence class; this can be based on a notion of minimal representation: within the IM class the canonical process is the process with the smallest number of states; within the OJC class, this is the process of lowest dimension. For instance for the illness-death process, the canonical IM is the three-state process $X$ (rather than the four-state process $W$) and the canonical OJC process is the bivariate ``basic'' counting process (rather than the three-dimensional process obtained form the Borgan representation).

In the next section we will derive the likelihood for the IM-OJC class of processes from increasingly complex schemes of observation and show in some examples that the rigorously derived likelihoods are the same as the heuristic ones.

\section{Derivation of the likelihood from
Jacod's formula}
\subsection{Likelihood for counting processes observed in continuous time}

From now on we adopt a somewhat more rigorous probabilistic formalism. We assume a probability space $(\Omega, \FF, P)$ is given and we consider a counting process $N$ defined on this space.
 Jacod (1975) has given the likelihood
ratio for observation of a counting process $N$ on $[0,C]$, that
is relative to the $\sigma$-field $\FF_C=\FF_0 \vee \N_C$ where
$\N_C=\sigma(N_{ju},0\le u\le C; j=1,\ldots,p)$. Aalen (1978),
based on results of Jacod \& Memin (1976), gave a simple form of
the likelihood ratio in the case of absolutely continuous
compensators by taking a reference probability  under which the
$N_j$'s are independent Poisson processes with intensity $1$. Here
we consider a multivariate counting process $N$ with components
$N_j$ which are $0-1$ counting processes. For such processes
 it is  more attractive to take a reference
probability $P_0$ under which the $N_j$'s are independent with
intensities $\lambda_j^0 (t)=1_{\{N_{jt-}=0\}}$; equivalently the
$T_j$'s are independent with exponential distributions with unit
parameter. The likelihood ratio for a probability $P_{\theta}$ (with $P_{\theta}<<P_0$) relative to $P_0$
is:
\begin{equation} \label{Jacod1}
\LL^{\theta}_{\FF_C}=\LL^{\theta}_{\FF_0}
\Bigl [\prod_{r=1}^{N_{.C}}\lambda^{\theta} _{J_r}(T_{(r)})\Bigr ]\exp
[-\Lambda^{\theta}_{.}(C) ]\prod_{j=1}^p e^{T_j\wedge C},
\end{equation}
where for each $r \in \{1 , \ldots ,N_{.C}\}$, $J_r$ is the unique
$j$ such that $\Delta N_{jT_{(r)}}=1$;
$N_{.t}=\sum_{j=1}^pN_{jt}$,
$\Lambda^{\theta}_{.}(t)=\sum_{j=1}^p\Lambda^{\theta}_{j}(t)$,
$\Lambda^{\theta}_{j}(t)=\int_0^t \lambda_j^{\theta}(u)du$ and $\lambda_j^{\theta}(t)$ is the intensity of $N_j$ under $P_{\theta}$. This
formula allows us to compute the likelihood for any multi-state model
once we have written it as a multivariate counting process. Within the OJC class, we denote $T_j$ the jump time of $N_j$; the
likelihood can then be written (in a more manageable form for applications)
in term of $\tilde T_j=\min (T_j, C)$ and $\delta_j=1_{\{T_j\le
C\}}$ as:
\begin{equation} \LL^{\theta}_{\FF_C}=\LL^{\theta}_{\FF_0} \Bigl
[\prod_{j=1}^{p}\lambda^{\theta}_{j}(\tilde T_j)^{\delta_{j}}\Bigr
]\exp [-\Lambda^{\theta}_{.}(C) ]\prod_{j=1}^p e^{\tilde T_j}.\end{equation}

  The term $\LL^{\theta}_{\FF_0}$ appears when there is some
information about the process of interest in the initial
$\sigma$-field at time 0; from now on we suppose that it is not
the case and so $\LL^{\theta}_{\FF_0}=1$ and this term disappears.
We may have to compute conditional likelihoods, in particular
likelihoods conditional on $\FF_{v_0}$ as we have seen it in
section 2. The conditional likelihood is simply
$\LL^{\theta}_{\FF_C|\FF_{v_0}}=\LL^{\theta}_{\FF_C}/\LL^{\theta}_{\FF_{v_0}}$.
In particular the likelihood $\LL^{\theta}_{\FF_{v_0}}$ on
$\{N_{v_0}=0\}$ is equal to $\exp [-\Lambda^{\theta}_{.}(v_0)
+pv_0]$ so that $\LL^{\theta}_{\FF_C|N_{v_0}=0}=\Bigl
[\prod_{j=1}^{p}\lambda^{\theta}_{j}(\tilde T_j)^{\delta_{j}}\Bigr
] \exp [-(\Lambda^{\theta}_{.}(C) -\Lambda^{\theta}_{.}(v_0)
)]\prod_{j=1}^p e^{\tilde T_j-v_0}$, (in this formula the the $\tilde T_j$ are the jump times possibily observed after $v_0$).

It is interesting to make the link with the heuristic likelihood
expressed in term of transition probabilities and intensities.
Dropping the multiplicative factor $\prod_{j=1}^p e^{\tilde T_j-v_0}$
(which does not depend on $\theta$) and rearranging the product,
the likelihood (\ref{Jacod1}) can be written:
$$\LL^{\theta}_{\FF_C}=\Bigl [\prod_{r=1}^{N_{.C}}\exp
-[\Lambda_{.}(T_{(r)})-\Lambda_{.}(T_{(r-1)})]\lambda
_{J_r}(T_{(r)})\Bigr
]\exp-[\Lambda_{.}(C)-\Lambda_{.}(T_{N_{.C}})],$$ still with the
convention $T_{(0)}=v_0$. We note that the number of jumps called $M$
in section 2.2 is precisely $N_{.C}$; on $\{X_{T_{(r)}}=x_r,
r=1,\ldots,m\}$ we have $\lambda
_{J_r}(T_{(r)})=\alpha_{x_{r-1},x_r}(T_{(r)})$ and $\exp
-[\Lambda_{.}(T_{(r)})-\Lambda_{.}(T_{r-1})]=\exp -[\int
_{T_{r-1}}^{T_{r}}\sum_{j}\lambda_j(u) du]= \exp -[\int
_{T_{r-1}}^{T_{r}}\sum_{h\ne x_{r-1}}\alpha_{x_{r-1}h}(u) du] $.
 The latter term is equal to $p_{\scriptscriptstyle x_{r-1},x_{r-1}}{\scriptstyle
(T_{(r-1)},T_{(r)}-)}$
 so that we
retrieve the expression given in 2.2.

 It
must be noted that in general the $\lambda^{\theta}_{j}(t)$'s and
$\Lambda^{\theta}_{j}(t)$'s may depend on what has been observed before
$t$; to be more explicit we can write $\lambda^{\theta}_{j}(t; \tilde T_l\wedge t,
l=1,\ldots p)$ and $\Lambda^{\theta}_{j}(C; \tilde T_l, l=1,\ldots
p)$.

The likelihood can be written as:
\begin{equation} \label{Jacod}\LL^{\theta}_{\FF_C}=f^{\theta}_C(\tilde T_1,\ldots,\tilde T_p)\prod_{j=1}^p e^{\tilde T_j}.\end{equation}
where
  $$f^{\theta}_C(s_1,\ldots,s_p)=\prod_{j=1}^{p}[\lambda^{\theta}_{j}(s_j;s_l\wedge s_j, l=1\ldots,p)]^{1_{\{s_j< C\}}}
\exp [-\Lambda^{\theta}_{.}(C;s_l \wedge C, l=1\ldots,p) ],$$
is the part of the likelihood which depends on $\theta$.
Note that in this expression of the likelihood we have get rid of the $\delta_j$'s, thus simplifying the notation for the developments of the next sections; in this expression it is considered that if $\tilde T_j=C$ the observation is right-censored; the case $\{T_j=C\}$ does not make problem because this event has probability zero with our assumptions and the likelihood is defined almost everywhere.
\subsection{Likelihood for the mixed discrete-continuous observation scheme: the case when
only one component is partially observed }
In this section, we
shall compute the likelihood for the following scheme of
observation of $N$: $N_1$ is observed at discrete times $v_0,
\ldots, v_m$ while $N_j, j=2, \dots, p$ are observed in continuous
time on $[0,C]$. The observation in this scheme is represented by
the $\sigma$-field $\Ob=\sigma(N_{1v_l}, l=0,\ldots ,m)\vee \G_C$
where $\G_C=\sigma(N_{ju}, j=2,\ldots, p; 0\le u \le C)$. We
obviously have $\Ob \subset \FF_C$ so that the observed likelihood
can be expressed as: $\LL^{\theta}_{\Ob}=\Ez [\LL^{\theta}_{\FF_C}|\Ob]$ where $\Ez$
means that we take expectation under $P_0$ defined in the previous
subsection. We note from formula (\ref{Jacod}) that
$\LL^{\theta}_{\FF_C}=f^{\theta}_C(\tilde T_1, \Gamma)g(\Gamma)e^{\tilde T_1}$, where $\Gamma= (\tilde T_2,\ldots,\tilde T_p)$, $g(\Gamma)= \prod_{j=2}^p e^{\tilde T_j}$ and $f^{\theta}_C(\tilde T_1, \Gamma)$ is a shortcut for $f^{\theta}_C(\tilde T_1, \tilde T_2,\ldots,\tilde T_p )$; note that $\Gamma$
is a $\G_C$-measurable random variable. From independence between
the $T_i$'s under $P_0$, we have independence between
$\sigma(N_{1v_l}, l=1,\ldots ,m)$ and $\G_C$ and  the computation
of the conditional expectation can be done using the
disintegration theorem (Kallenberg, 2001) which yields:
$$\LL^{\theta}_{\Ob}=g(\Gamma)\int f^{\theta}_C(s,\Gamma)e^{s\wedge
C} \nu(ds),$$ where $\nu(.)$ is a regular version of the
law of $T_1$ given $\Ob$ which, by independence equals a version
of the law of $T_1$ given $\sigma(N_{1v_l}, l=0,\ldots ,m)$. We
decompose $\LL^{\theta}_{\Ob}$ on atoms (see a definition in section 4.4) of $\sigma(N_{1v_l}, l=0,\ldots
,m)$ as :
\begin{equation} \label{decomp} \LL^{\theta}_{\Ob}=\sum_{l=1}^m
1_{\{v_{l-1}<T_1\le v_l\}}\Ez [\LL^{\theta}_{\FF_C}|\Ob] +1_{\{T_1 >
v_m\}}\Ez [\LL^{\theta}_{\FF_C}|\Ob],\end{equation} with the convention $v_{0}=0$. Using the fact that
$\nu
(ds)=\frac{1_{(v_{l-1},v_l]}(s)e^{-s}}{(e^{-v_{l-1}}-e^{-v_l})}ds$
gives the law of  $T_1$ given $(T_1 \in (v_{l-1},v_l])$, we obtain
for the first terms:
\begin{eqnarray}1_{\{v_{l-1}<T_1\le v_l\}}\Ez
[\LL^{\theta}_{\FF_C}|\Ob] & = & 1_{\{v_{l-1}<T_1\le v_l\}}g(\Gamma)\Ez
[f^{\theta}_C(T_1,\Gamma)|T_1 \in ]v_{l-1}; v_l]]\nonumber\\ & = &
1_{\{v_{l-1}<T_1\le
v_l\}}\frac{g(\Gamma)}{e^{-v_{l-1}}-e^{-v_l}}\int_{v_{l-1}}^{v_l}
f^{\theta}_C(s,\Gamma)ds\label{vi}.\end{eqnarray}
\noindent We consider now the last term of (\ref{decomp}). $\nu (ds)=1_{(v_m,+ \infty
]}(s)e^{-s}e^{v_{m}}ds$ is a regular version of the law of $T_1$
given $(T_1>v_m)$ and we can write
\begin{eqnarray}1_{\{T_1 >
v_m\}}\Ez[\LL^{\theta}_{\FF_C}|\Ob] & = & 1_{\{T_1 > v_m\}}g(\Gamma)\Ez
[f^{\theta}_C(T_1,\Gamma)|T_1 > v_m]\nonumber\\ & = &
1_{\{v_{m}<T_1\}}g(\Gamma)e^{v_{m}}\left[\int_{v_{m}}^{C}
f^{\theta}_C(s,\Gamma)ds + \int_{C}^{+\infty}
f^{\theta}_C(C,\Gamma)e^{C-s}ds\right]\nonumber\\ & = & 1_{\{T_1 >
v_m\}}g(\Gamma)e^{v_m}\left[\int_{v_{m}}^{C}
f^{\theta}_C(s,\Gamma)ds+f^{\theta}_C(C,\Gamma)\right]\label{vm}.
\end{eqnarray}
Combining equations (\ref{decomp}), (\ref{vi}) and
(\ref{vm}) we have proved

\begin{Lemma}\label{pas1} For $N_1$
observed at discrete times $v_0, \ldots,v_m$ and $N_j$, ${(2\leq j
\leq p)}$ observed in continuous time the likelihood is given by
\begin{eqnarray} \LL^{\theta}_{\Ob}& = & \sum_{l=1}^m 1_{\{v_{l-1}<T_1\le
v_l\}}\frac{g(\Gamma)}{e^{-v_{l-1}}-e^{-v_l}}\int_{v_{l-1}}^{v_l}
f^{\theta}_C(s,\Gamma)ds  \label{decomplete}\\ &   &  +1_{\{T_1 >
v_m\}}g(\Gamma)e^{v_m}\left[\int_{v_{m}}^{C}
f^{\theta}_C(s,\Gamma)ds+f^{\theta}_C(C,\Gamma)\right]   \nonumber
\end{eqnarray} \end{Lemma}

\noindent
 {\em Example: The Markov illness-death
model.} Let us apply the formula to a Markov illness-death model;
for brevity we shall compute the likelihood only on the event $
\{v_{l-1}<T_1\le v_l\} \cap \{T_2 > C\}$ (that is, in
epidemiological language, when a subject has been first seen ill
at $v_l$ and was still alive at $C$). On
this event ${N_{.C}=1}$ and in this case:
$f^{\theta}_C(T_1,\Gamma)=\lambda_1(T_1)e^{
-\Lambda^{\theta}_{.}(C)}$. (here $\Gamma=\tilde T_2$ and $\tilde T_2
=C$ on this event). We decompose the cumulative total
intensity as: ${ \Lambda^{\theta}_{.}(C)}=
\Lambda^{\theta}_{.}(v_{l-1})- \Lambda^{\theta}_{.}(v_{l-1})+\Lambda^{\theta}_{.}(T_1)- \Lambda^{\theta}_{.}(T_1)+
\Lambda^{\theta}_{.}(v_l)- \Lambda^{\theta}_{.}(v_l)+
\Lambda^{\theta}_{.}(C)$. We have in the Markov model (from
formulae (\ref{int1}) and (\ref{int2})) that
$\Lambda^{\theta}_{.}(v_{l-1})=A_{0.}(0,v_{l-1})$, $\Lambda^{\theta}_{.}(T_1)-\Lambda^{\theta}_{.}(v_{l-1})=A_{0.}(v_{l-1},T_1)$, $\Lambda^{\theta}_{.}(v_l)-\Lambda^{\theta}_{.}(T_1)=A_{12}(T_1,v_l)$ and $\Lambda^{\theta}_{.}(C)-\Lambda^{\theta}_{.}(v_l)=A_{12}(v_l,C)$, where $A_{ij}(a,b)=\int _a^b \alpha_{ij}(s) ds$. So that, taking out of the
integral the terms which do not depend on $s$ we obtain:
 $$\int_{v_{l-1}}^{v_l} f^{\theta}_C(s,\Gamma)ds=e^{
-A_{0.}(0,v_{l-1})}\Bigl [\int_{v_{l-1}}^{v_l}e^{
-A_{0.}(v_{l-1},s)} \alpha_{01}(s)e^{
-A_{12}(s,v_l)}ds \Bigr ]e^{ -A_{12}(v_l,C)}.$$
This is  a product of three terms in which we recognize
$p_{00}(0,v_{l-1})$, $p_{01}(v_{l-1},v_l)$ and $p_{11}(v_l,\tilde
T)$ (Commenges {\em et al.}, 2004). Conditioning on $\{N_{v_0}=0\}$ we have to divide by $e^{-A_{0.}(0,v_{0})}$ and  we retrieve formula (\ref{liki}).

\subsection{Likelihood for the mixed discrete-continuous
observation scheme: the case when several components are partially
observed } In this section, we shall first compute the likelihood
for the following scheme of observation of $N$: $N_1$ is observed
at discrete times $v^1_0, \ldots, v^1_{m_1}$,  and $N_2$ is
observed at discrete times $v^2_0, \ldots, v^2_{m_2}$ while $N_j,
j=3, \dots, p$ are observed in continuous time on $[0,C]$ (see section 4.5 for the problem raised by the randomness of the number of observations). The
observation in this scheme is represented by the $\sigma$-field
$\Ob=\sigma(N_{1v^1_{l_1}}, N_{2v^2_{l_2}} 0 \leq l_1\leq m_1 ,0
\leq l_2\leq m_2)\vee \G^2_C$ where $\G^2_C=\sigma(N_{ju},
j=3,\ldots, p; 0\le u \le C)$. Using (\ref{Jacod}) we write
$\LL^{\theta}_{\FF_C}=f^{\theta}_C(\tilde T_1,\tilde T_2,\Gamma^2)g(\Gamma^2)e^{\tilde T_1+\tilde T_2}$,
 where
$\Gamma^2$ is a $\G^2_C$-measurable random variable. If we note
$\Ob'=\sigma(N_{1v^1_{l_1}}, l_1=0,\ldots ,m_1)\vee \G_C$ (which
we denoted $\Ob$ in the previous subsection), we obviously have
$\Ob \subset \Ob' \subset \FF_C$ so that the observed likelihood
can be expressed as $\LL^{\theta}_{\Ob}=\Ez[\Ez[\LL^{\theta}_{\FF_C}|\Ob']|\Ob]$.
 Lemma \ref{pas1} gives an expression of
$\Ez[\LL^{\theta}_{\FF_C}|\Ob']$ and using the same properties as above and
the convention  $v^1_0=0$ and $v^2_0=0$, we can write :
\begin{eqnarray}
& \LL^{\theta}_{\Ob}&  \nonumber \\ &  = & \sum_{l_1=1}^{m_1}
\sum_{l_2=1}^{m_2} \frac{ 1_{\{v^1_{l_1-1}<T_1\le v^1_{l_1}\}}}
{e^{-v^1_{l_1-1}}-e^{-v^1_{l_1}}}1_{\{v^2_{l_2-1}<T_2\le
v^2_{l_2}\}}g(\Gamma^2) \Ez \left[\int_{v^1_{l_1-1}}^{v^1_{l_1}}
e^{\tilde T_2} f^{\theta}_C(s,\tilde T_2 ,\Gamma^2)ds
|\Ob\right] \label{ij}  \nonumber\\
 &  + &   \sum_{l_2=1}^{m_2}
1_{\{v^2_{l_2-1}<T_2\le v^2_{l_2}\}}1_{\{T_1 >
v^1_{m_1}\}}g(\Gamma^2)e^{v^1_{m_1}}\times \nonumber\\ &    &
\;\;\;\;\; \Ez \left[e^{\tilde T_2} [\int_{v^1_{{m_1}}}^{C}
f^{\theta}_C(s_1,\tilde T_2,\Gamma^2)ds_1+f^{\theta}_C(C,\tilde T_2,\Gamma^2)] |\Ob\right] \nonumber\\
&  + &   \sum_{l_1=1}^{m_1} \frac{1_{\{v^1_{l_1-1}<T_1\le
v^1_{l_1}\}}}{e^{-v^1_{l_1-1}}-e^{-v^1_{l_1}}}1_{\{T_2 >
v^2_{m_2}\}}g(\Gamma^2)e^{v^2_{m_2}} \Ez
\left[\int_{v^1_{l_1-1}}^{v^1_{l_1}} e^{\tilde T_2}
f^{\theta}_C(s_1,\tilde T_2 ,\Gamma^2)ds_1 |\Ob\right]
\nonumber \\
& + & 1_{\{T_1 > v^1_{m_1}\}}1_{\{T_2 >
v^2_{m_2}\}}e^1{v_{m_1}}\Ez \left[e^{\tilde T_2}
[\int_{v^1_{{m_1}}}^{C} f^{\theta}_C(s_1,\tilde T_2,\Gamma^2)ds_1
+f^{\theta}_C(C,\tilde T_2,\Gamma^2)] |\Ob\right].
\nonumber \end{eqnarray}
We compute each term of the right hand of
this equality using the law and the independence of the $T_i$'s
under $P_0$ as in the previous subsection. Dropping the terms
which not depend on $\theta$ (as $g(\Gamma^2)$) we have the
\begin{Lemma}\label{pas2} For $N_1$ observed
at discrete times $v^1_0, \ldots, v^1_{m_1}$,  and $N_2$ observed
at discrete times $v^2_0, \ldots, v^2_{m_2}$ and $N_j$, ${(3\leq j
\leq p)}$ observed in continuous time the likelihood is given by
\begin{eqnarray*} & \LL^{\theta}_{\Ob}&  \nonumber \\ &  = &
\sum_{l_1=1}^{m_1} \sum_{l_2=1}^{m_2} \frac{
1_{\{v^1_{l_1-1}<T_1\le v^1_{l_1}\}}}
{(e^{-v^1_{l_1-1}}-e^{-v^1_{l_1}})}\frac{1_{\{v^2_{l_2-1}<T_2\le
v^2_{l_2}\}}}{(e^{-v^2_{l_2-1}}-e^{-v^2_{l_2}})}
\int_{v^1_{l_1-1}}^{v^1_{l_1}}\int_{v^2_{l_2-1}}^{v^2_{l_2}}
f^{\theta}_C(s_1,s_2,\Gamma^2)ds_1ds_2 \nonumber \\
&  + & \sum_{l_2=1}^{m_2} e^{v^1_{m_1}}1_{\{T_1 >
v^1_{m_1}\}}\frac{1_{\{v^2_{l_2-1}<T_2\le
v^2_{l_2}\}}}{e^{-v^2_{l_2-1}}-e^{-v^2_{l_2}}}
\int_{v^2_{l_2-1}}^{v^2_{l_2}}(\int_{v^1_{m_1}}^{C}
f^{\theta}_C(s_1,s_2,\Gamma^2)ds_1+f^{\theta}_C(C,s_2,\Gamma^2) )ds_2
\nonumber\\
&  + &  \sum_{l_1=1}^{m_1}  \frac{ 1_{\{v^1_{l_1-1}<T_1\le
v^1_{l_1}\}}}
{e^{-v^1_{l_1-1}}-e^{-v^1_{l_1}}}e^{v^2_{m_2}}1_{\{T_2
> v^2_{m_2}\}} \int_{v^2_{l_2-1}}^{v^2_{l_2}}(\int_{v^1_{m_1}}^{C}
f^{\theta}_C(s_1,s_2,\Gamma^2)ds_2+f^{\theta}_C(s_1,C,\Gamma^2) )ds_1
\nonumber\\ & + & 1_{\{T_1 > v^1_{m_1}\}}1_{\{T_2 >
v^2_{m_2}\}}e^{(v^1_{m_1}+v^2_{m_2})} \times
\left[\int_{v^2_{m_2}}^{C} \int_{v^1_{m_1}}^{C}
f^{\theta}_C(s_1,s_2,\Gamma^2)ds_1ds_2 \right.\nonumber
\\ &   & \;\;\;\; + \left.
\int_{v^1_{m_1}}^{C}f^{\theta}_C(s_1,C,\Gamma^2)ds_1+\int_{v^2_{m_2}}^{C}
f^{\theta}_C(C,s_2,\Gamma^2)ds_2 +
f^{\theta}_C(C,C,\Gamma^2)\right].\nonumber
\end{eqnarray*} \end{Lemma}
We can prove by induction a formula for the
likelihood when $k$ components of the processes are partially
observed. The high number of possibilities for the realization of
the $(T_i)_{1 \le i \le k }$ in the intervals
$[v^i_{l_i};v^i_{l_i-1}]$ or $[v^i_{m_i};+\infty[$ makes the
formula complicated. Lemma \label{local} in the next sub-section makes it unnecessary to give such a cumbersome formula.

\subsection{General formula of the likelihood via its local representation}
Here we develop an approach
which is closer to the statistician's point of view and which is
more general. We begin with a Lemma.
\begin{Lemma}\label{local}
Consider two observation schemes of $N$ yielding observed
$\sigma$-fields $\Ob$ and $\tilde \Ob$; consider an event $A$ such
that $P_0(A)>0$, $A\in \Ob \cap \tilde \Ob$ and $A\cap \Ob = A
\cap \tilde \Ob$; then $\LL^{\theta}_{\Ob}=\LL^{\theta}_{\tilde \Ob}$ a.s. on $A$.
\end{Lemma}

\noindent {\em Proof.} Remember that $\LL^{\theta}_{\Ob}=\Ez[\LL^{\theta}_{\FF}|\Ob]$ and
$\LL^{\theta}_{\tilde \Ob}=\Ez[\LL^{\theta}_{\FF}|\tilde \Ob]$. Direct application
of a lemma of local equality of conditional expectations
(Kallenberg, 2001: Lemma 2, Chapter 6) gives the result.

For instance consider the case where $N_1$ is observed at discrete
times $v_0,\ldots,v_m$ and $N_j, j>1$ are observed in continuous time as in section 4.1; take
$A=\{T_1\in (v_{l-1},v_l]\}$.We have
$\LL^{\theta}_{\Ob}=\frac{g(\Gamma)}{e^{-v_{l-1}}-e^{-v_l}}\int_{v_{l-1}}^{v_l}f^{\theta}(s,\Gamma)ds$,
a.s. on $A$ (this can be seen by multiplying both sides of
equation (\ref{decomplete}) by $1_A$ which is equal to $1$ on $A$).
From a statistician's point of view we can say that if  $A$ happens
 the likelihood takes that particular form which does not depend on the other values of
the observation times ($v_{l'}, l'\ne l, l-1$). So it is obvious
either directly, or using Lemma \ref{local} that if we consider
another scheme of observation where $N_1$ is observed at
$w_1=v_{l-1}$ and $w_2=v_{l}$ leading to the likelihood
$\LL^{\theta}_{\tilde \Ob}$, we have $\LL^{\theta}_{\Ob}=\LL^{\theta}_{\tilde \Ob}$ a.s. on
$A$.

\vspace{0.25cm} \noindent
This leads to an extension of the field
of application of the formulae of the likelihood for incomplete
data. A more general scheme of observation
 can be considered:

\begin{Definition}[The deterministic GCMP]
A deterministic GCMP is a scheme of observation for a multivariate process $X=(X_1,\ldots,X_p)$ specified by a response function $r(.)=(r_1(.),\ldots,r_p(.))$, where the $r_j(.)$'s take values $0$ or $1$,  such that $X_{jt}$ is observed at time $t$ if and only if $r_{j}(t)=1$, for $j=1,\ldots,p$.
\end{Definition}
This general (deterministic) scheme applied to $N$ allows each
component $N_j$ to be observed in continuous time over some
windows and in discrete time over other windows. Within this
scheme (assuming a family of equivalent probability measures), the
likelihood is given by $\Ez[\LL^{\theta}_{\FF}|\Ob]$ where $\Ob
=\sigma(r_{j}(t)N_{jt} ;0 \leq t \leq C; j \in \{1, \ldots ,p\})$.
Lemma \ref{local} will help us to give a simple expression of this
likelihood if we can find a finite class of events $(A_k)$ which form a partition of $\Omega$ and on which the likelihood is relatively easy to compute. If all the observed events had a positive (and bounded away from zero) probability, the class of atoms of $\Ob$ would yield a natural finite partition of $\Omega$ (a P-atom of a $\sigma$-field is a set $A$ belonging to it, such that $P(A)>0$ and if $B\subset A$, then $P(B)=0$ or $P(A\cap \bar B)=0$). For instance this would be the case if all the components of $N$ were observed in discrete time. As soon as one component may be observed in continuous time, there may still exist atoms, but we do not have a partition of $\Omega$ with atoms. This leads us to define a class of pseudo-atoms which is a partition of $\Omega$. For $0 \le j \le p$, we denote :

\begin {itemize}
\item[-]$v^j_0=0$ and $\{v^j_1,\ldots,v^j_{m_j} \}$ the finite set of times of
discontinuities of the function $r_{j}(t)$ on $[0,C]$ such
that $\{v^j_0=0 <v^j_1 <\ldots<v^j_{m_j} \}$.

\item[-] if $m_j>0$ we define $D_{jl}=\{T_j \in (v^j_{l-1}, v^j_l]\} \cap \{r_{ju}=0, u \in
(v^j_{l-1},v^j_l)\}$, $l \geq 1$. For some $l$, $D_{jl}$ is empty
otherwise it is an atom of $\Ob_j=\sigma(r_{j}(t)N_{jt} ;0 \leq t
\leq C)$.
\item[-]$E_{j}=(\{T_j>v^j_{m^j}\}\cap \{r_{ju}=0, u>v^j_{m_j}\}) \cup \{T_j>C\}$ which also is either empty or an atom of $\Ob_j$.

\item[-] $C_j=\{T_j\in \stackrel{\circ}{\widehat{r_j^{-1}(1)}}\}$ where
$\stackrel{\circ}{\widehat{r_j^{-1}(1)}}$ denotes the interior (in
topological sense) of the $\Re$-subset $r_j^{-1}(1)$ in
$[0,v^j_{m_j}]$. On $C_j$, the event $T_j$ is exactly observed and
$C_j$ is the complementary set of $\left[\left(\cup_{1 \le l \le =
\infty} D_{jl}\right) \cup E_j\right]$ in $\Omega$.
\end{itemize}
 $\Omega$ can be partitioned into the disjoint sets $D_{jl},
l \geq 1$,  $C_j$ and $E_j$ for each $j$
 (there is
no such $D_{jl}$ if $N_j$ is observed in continuous time over
$[0,C]$; also $C_j$ is empty if $N_j$ is only observed at
discrete times). A finer partition can be obtained by the intersection of these $p$ partitions.

\begin{Definition}[Pseudo-atoms]
In the deterministic GCMP framework, we call pseudo-atoms of
$\Ob$, a set $A=\cap _{j=1}^p A_j$, where $A_j=D_{jl}$ for some
$l$,  $A_j=C_{j}$ or $A_j=E_j$ and such that
$P_0(A)>0$.
\end{Definition}
\noindent {\em Remark.} If $A$ is a pseudo-atom of $\Ob$ we have
$A\in \Ob$. The class of pseudo-atoms $(A_k)$ form a partition of
$\Omega$: $\Omega= \cup_k A_k; A_k\cap A_{k'}=\emptyset$ if $ k
\ne k'$.

\noindent{\em Example: The Illness-death model with hybrid observation scheme for illness: pseudo-atoms.}
 Consider an illness-death model in which illness represents a complication of a disease. The occurrence of
  the complication (illness) is observed in continuous time during a sojourn in hospital from
  time $0$ to time $v_1$. After the
  hospitalization,  the complication is diagnosed at planned visits
  at times $\{v_2, v_3, \ldots,v_m \leq C\}$. We suppose that if
  death occurs during the study (i.e. between $[0,C]$) its
  time can be retrieved exactly.
  So the process of observation is given
  by $ \left\{ \begin{array}{rcl}
          r_1(t)& = & 1_{[0,v_1)}(t) + \sum_{i=2}^m 1_{\{t=v_i\}} \\
         r_{2}(t)  & = & 1 \ \  0 \leq t \leq C
       \end{array}\right.$
With this scheme, the times of discontinuities and the sets defined previously are :

\begin{itemize}
\item[-]For the first component $r_1(.)$, $m_1=m$ and $v^1_k=v_k$.
For the second one, $m_2=0$.

\item[-] $D_{11}$ is empty and for $2 \leq l \leq m$, $D_{1l}= \{T_1 \in (v_{l-1}, v_l] \}$.
 If the event $D_{1l}$ occurs for some $l \geq 2$ it means that
the complication is diagnosed for the first time at the visit $v_l$.
 Since $m_2=0$ there is no  $D_{2l}$.

\item[-]$E_{1}=\{T_1 >v_{m}\}$ and $E_2= \{T_2>C\}$
\item[-] $C_1=\{T_1\le v_1\}$. If $C_1$ occurs it means
that the illness is diagnosed during the sojourn in hospital.
 We have $C_2=\{T_2  \leq C\}$.
\end{itemize}

Note that $\Omega=C_1 \cup (\cup_{l=2}^m D_{1l}) \cup E_{1}$ gives
a partition of $\Omega$ and $\Omega= C_2 \cup E_2$ gives an other
one and the pseudo-atoms are the intersection of sets of these two partitions of positive probability.
 Now we can give the likelihood by describing it on any pseudo-atom.

\begin{Theorem}\label{OA}
Consider the deterministic GCMP scheme specified by $r_j(.),
j=1,\ldots, p$, and a pseudo-atom of $\Ob$, $A=\cap _{l\in L}\{T_l
\in (v^l_1,v^l_2]\}\cap _{l\in L'}\{T_{l'}>v^{l'}_{m_{l'}} \}\cap
_{l\in \bar L \cap \bar L'} \{T_l\in r_l^{-1}(1)\}$ where $L$ and
$L'$ are disjoint subsets of $\{1,\ldots p\}$. On $A$ the
likelihood $\LL^{\theta}_{\Ob}$ is equal to that of a MDCO scheme where
$N_l$ is observed at times $v_1^l$ and $v_2^l$ for $l\in L$, at
$v^{l'}_{m_{l'}}$, for $l'\in L'$ and in continuous time for
$l''\in L''=\bar L\cap \bar L'$. Without loss of generality, we
can suppose that $L''=\{k+1,\ldots,p\}$. On A this likelihood is
equal to:

\begin{eqnarray*}
 & \LL^{\theta}_{\Ob} & = \frac{e^{(\sum_{l' \in L'}{v^{l'}_{m_{l'}}})}}
 {\prod_{l\in L} (e^{-v^l_{1}}- e^{-v^l_{2}})} \left[\oint
\bigg [\prod_{l\in L} 1_{{\{v^l_{1}<s_l \leq v^l_{2}\}}}\prod_{l'\in
L'} 1_{\{v^{l'}_{m_{l'}}<s_{l'}\leq C
\}}\Bigr ]f_C^{\theta}(s_1,\ldots,s_k,\Gamma_k)
\prod_{l \in L \cup L'}ds_l\right.\label{1} \\
&  + &   \sum_{i\in L'} \oint \bigg [\prod_{l\in L} 1_{\{v^l_{1}<s_l
\leq v^l_{2}\}} \prod_{l'\in L', l'\neq i}
1_{\{v^{l'}_{m_{l'}}<s_l \leq C
\}}\Bigr ]f_C^{\theta}(s_1,\ldots,\underbrace{C}_i,\ldots,s_k,\Gamma_k)
\prod_{l\neq i }ds_l\nonumber \\
 &  + &   \sum_{(i_1 \neq i_2) \in L'^2} \oint \bigg [\prod_{l\in L}
1_{\{v^l_{j_l-1}<s_l \leq v^l_{j_l}\}}\prod_{l'\in L', l'\neq
i_j}
1_{\{v^{l'}_{m_{l'}}<s_l\}}\Bigr ]f_C^{\theta}(s_1,.,\underbrace{C}_{i_1},.,\underbrace{C}_{i_2}.,s_k,\Gamma_k)
\prod_{l\neq i_j }ds_l\nonumber \\ \\ &  + & \left.  \ldots \;\
\right]
\end{eqnarray*}
\end{Theorem}

\noindent {\em Proof.}  Let $\tilde \Ob$ the observed
$\sigma$-field in the MDCO scheme of the Theorem. We have $A\in
\Ob$, $A\in \tilde \Ob$ and $A\cap \Ob=A\cap \tilde \Ob$. Lemma
\ref{local} gives us $\LL^{\theta}_{\Ob}=\LL^{\theta}_{\tilde \Ob}$ a.s on $A$. The
value of the likelihood on this event can be derived using the
technique of the preceding sub-section.

\noindent{\em Example: The Illness-death model with hybrid observation scheme for illness: likelihood.}
The partition of $\Omega$ given by the class of pseudo-atoms is $\{C_1
\cap C_2 ;  D_{1l} \cap C_2 , l=2,\ldots,m; E_{1} \cap C_2 ; C_1 \cap E_2 ;
D_{1l} \cap E_2, l=2,\ldots,m; E_{1} \cap E_2\}$. In this case we have
$f^{\theta}_C(s_1,s_2)=\lambda^{\theta}_1(s_1;s_1, s_2\wedge s_1)^{1_{\{s_1 < C\}}}\lambda^{\theta}_2(s_2;s_1\wedge s_2, s_2)^{1_{\{s_2 < C\}}}  exp[-\Lambda^{\theta}_{.}(C;s_1\wedge C ,s_2\wedge C )].$
The likelihood is easy to write on each
of the pseudo-atoms.
\begin{itemize}
\item[-] On $\{C_1\cap C_2\}$ , all the components are observed in continuous time ($L$ and $L'$ are
empty) : $\LL^{\theta}_{\Ob}=\LL^{\theta}_{\FF_C}$;
\item[-] On $D_{1l} \cap C_2$, $L=\{1\}$ and $L''=\{2\}$ ,
$\LL^{\theta}_{\Ob}=\frac{1}{(e^{-v^1_{l-1}}-
e^{-v^1_{l}})}\int_{v^1_{l-1}}^{v^1_{l}}f_C^{\theta}(s_1,T_2)ds_1$;
\item[-] On $E_1 \cap C_2$, $\LL^{\theta}_{\Ob}=e^{v_m} \left[ \int_{v_m}^C f_C^{\theta}(s_1,T_2)ds_1 + f_C^{\theta}(C,T_2)
\right]$ We can remark that the intensity $\lambda_1^{\theta}(t)$
($\lambda_2^{\theta}(t)$) of occurrence of the illness (death)
vanishes after the death and then the likelihood is equal
to $\LL^{\theta}_{\Ob}=e^{v_m} \left[ \int_{v_m}^{T_2} f_
C^{\theta}(s_1,T_2)ds_1 + f_C^{\theta}(T_2,T_2) \right]$;
 \item[-] On  $C_1 \cap E_2$ all the components are observed in continuous
 time on $[0,C]$ : $\LL^{\theta}_{\Ob}=\LL^{\theta}_{\FF_C}$;
 \item[-] On $D_{1l} \cap E_2$, $L=\{1\}$ and $L''=\{2\}$,
$\LL^{\theta}_{\Ob}=\frac{1}{(e^{-v^1_{l-1}}-
e^{-v^1_{l}})}\int_{v^1_{l-1}}^{v^1_{l}}f_C^{\theta}(s_1,C)ds_1$;
\item[-] On $E_1 \cap E_2$,  one finds $\LL^{\theta}_{\Ob}=e^{(v_m+C)} \left[ \int_{v_m}^C f_C^{\theta}(s_1,C)ds_1 + f_C^{\theta}(C,C)
\right]$.
\end{itemize}

\subsection{Ignorability of the stochastic censoring by death}
In many examples of interest in epidemiology, one of the processes, say $N_p$ counts death, and the observation of $N_1, \ldots, N_{p-1}$ may be right censored by the time of death $T_p$ (in addition to other types of coarsening). For instance suppose it has been planned to visit a subject at discrete visit times
$v^1_0,\ldots,v^1_{m_1}\le C$ to observe the process $N_1$ (representing an illness status), the last visit time is necesssarily random: it is $v_{M_1}$ where $M_1=\max_l (l: v^1_l <\tilde T_p)$. So for treating this problem we have to resort to the stochastic GCMP: this is the same definition as the deterministic GCMP except that the response process is stochastic and will be denoted by $R$. It is assumed that $R$ is observed. In this problem we have to consider the $\sigma$-field $\FF=\R \vee \N_C$ and the observed $\sigma$-field is
$\Ob=\sigma(R_{jt}, R_{jt}N_{jt}, 0\le t\le C; j=1\ldots,p)$. Consider the case where the event $\{R=r\}$ has a positive probability. If the mechanism leading to incomplete data had been deterministically fixed to be equal to $r$, the observed $\sigma$-field would be $\tilde \Ob= \sigma(r_j(t)N_{jt}, 0\le t\le C; j=1\ldots,p)$. The mechanism leading to incomplete data is ignorabble if on the event  $\{R=r\}$ using $\LL^{\theta}_{\tilde \Ob}$ leads to the same inference as using $\LL^{\theta}_{ \Ob}$. Commenges \& G\'egout-Petit (2005) gave general conditions of ignorability.

 Let us treat here the specific problem in which we assume that the stochastic nature of $R$ comes only from the right censoring of the other processes by death. In that case the response processes can be written:
$$R_{jt}=r_j(t) 1_{\{N_{pt}=0\}}, 0\le t\le C; j=1,\ldots,p-1,$$
and we assume that $R_{pt}=1, 0\le t \le C$; the $r_j(t)$'s are deterministic functions as in Definition 2.
Define $M_j=\max_l (l: v^j_l <\tilde T_p)$ and consider  the event $A=\cap_{j=1}^{p-1} \{M_j=\tilde m_j\}$ for some choice of the $\tilde m_j$'s compatible with the deterministic discontinuities $v^j_l$'s; $A$ has then a positive probability if the intensity of $N_p$ does not vanish on $[0,C]$. Since the $R_{j}$ are functions of $N_p$ the observed $\sigma$-field can be written:
$\Ob=\sigma(r_j(t)1_{\{N_{pt}=0\}}N_{jt}, j=1\ldots,p-1; N_{pt}, 0\le t\le C)$.
 For  each $A$ we can define determistic functions as: $\tilde r_j(t)=r_j(t), t \le v_{\tilde m_j}; \tilde r_j(t)=r_j(t+), t >  v_{\tilde m_j}; j=1,\ldots,p-1$. For the deterministic GCMP specified by the $r_j(t)$'s the observed $\sigma$-field is $\tilde \Ob= \sigma(\tilde r_j(t)N_{jt}, 0\le t\le C; j=1\ldots,p)$, where $\tilde r_p(t)=1$, $0\le t\le C$. Because the stochastic right censoring depends only on $N_p$ which is observed, it is clear that $A\in \tilde \Ob$. Moreover because the intensities of the processes $N_j, j=1,\ldots, p-1$ are zero after death we have $\tilde \Ob=\Ob$ on $A$. Thus we can apply Lemma \ref{local} to find that on $A$  we have $\LL^{\theta}_{\tilde \Ob}=\LL^{\theta}_{ \Ob}$. The consequence is that in that case we can use the formulae derived in the deterministic framework (while it can not be deterministic), interpreting $v^j_{m_j}$ as the last discontinuity of $r_j(t)$ before death.

\subsection{Extension to general multi-state models}
Any multi-state process can be represented by an informationally equivalent counting process (obvious by the Borgan representation). The approach we have developed for deriving the likelihood in the GCMP context can be applied to general counting processes: that is compute the likelihood ratio for continuous-time observation and take the conditional expectation given $\Ob$. However it seems nearly impossible to obtain a general formula in that case, in particular because the number of transitions which may occur between two observation times is not bounded. In applications this is of course not realistic. For instance consider a two-state reversible model for diarrhoea: state $0$: no diarrhoea; state $1$ : diarrhoea. By definition a period of diarrhoea lasts a certain time, say one day. In that case the number of transitions between two observation times is bounded. It is then possible to cast the problem in the OJC framework because a counting process making at most $k$ jumps can be represented by a $k$-dimensional OJC process.

\section{Illustration on dementia-institution-death models}
\subsection{Multi-State and counting process models}
A five-state non-homogeneous Markov model
for dementia, institutionalization and death was proposed by
Commenges \& Joly (2004) (see Figure 1); we present this model
with somewhat different notations to be in agreement with our
general notations of section 4. The first aim was to estimate the
eight transition intensities. It was  proposed to make no
parametric assumption on
$\alpha_{01}(t),\alpha_{02}(t),\alpha_{04}(t)$ but to relate
parametrically the other transition intensities to these three
basic intensities. Proportionality of the transition intensities
toward dementia was assumed:
$\alpha_{23}(t)=\alpha_{01}(t)e^{\eta_2^1},$ as well as
proportionality of  the transition intensities toward
institution: $\alpha_{13}(t)=\alpha_{02}(t)e^{\eta_1^2}.$ It was
assumed that the transitions toward death are all proportional:
$\alpha_{14}(t)=\alpha_{04}(t)e^{\eta_1^3};$
$\alpha_{24}(t)=\alpha_{04}(t)e^{\eta_2^3};$
$\alpha_{34}(t)=\alpha_{04}(t)e^{\eta_{12}^3}.$
As shown in
section 3 this model can be represented by a trivariate process
$N=(N_1,N_2,N_3)$ where $N_1$ counts dementia, $N_2$ counts
institutionalization and $N_3$ counts death; the value of the
five-state process $X$ at $t$  can be represented in base 2 as:
$X_t=\min (N_{3t}N_{2t}N_{1t},4)$. The processes $X$ and $N$ are
equivalent in the sense that they generate the same filtration.
Moreover if under a probability measure $P$, $X$ is Markov with
transition intensities $(\alpha_{hj}, h=0,\ldots,4; j=0,\ldots,4)$
the intensity of $N$ can be deduced by the general formula
(\ref{intens-intens}). The model which was proposed by Commenges \& Joly
(2004) is a Markov semi-parametric multiplicative model; we write the intensities as functions of $T_j, j=1,2,3$ and we can verify that the intensities at $t$ only depend on $T_j\wedge t, j=1,2,3$ (because for instance $\{T_{j}\ge t\}=\{T_{j}\wedge t \ge t\}$); they are:
\begin{eqnarray} \label{Markov} \lambda^{\theta}_1(t;T_1,T_2,T_3)&=&\alpha_{01}(t)1_{\{T_{1}\ge t\}}1_{\{T_{3}\ge t\}}e^{\eta_2^1 1_{T_2<t}} \nonumber \\
\lambda^{\theta}_2(t;T_1,T_2,T_3)&=&\alpha_{02}(t)1_{\{T_{2}\ge t\}}1_{\{T_{3}\ge t\}}e^{\eta_1^2 1_{T_1<t}}\\
\lambda^{\theta}_3(t;T_1,T_2,T_3)&=&\alpha_{04}(t)1_{\{T_{3}\ge
t\}}e^{\eta_1^31_{T_1<t}+\eta_2^3 1_{T_2<t}+\eta^3_{12}
1_{T_1<t}1_{T_2<t}}\nonumber ,\end{eqnarray} where $\theta=(\eta ;
\alpha_{0j}(.), j=1,2,4)$ and $\eta$ represents the vector of
parameters named with this letter. It is attractive to consider some
non-Markovian models and, including an explanatory variable $Z_i$.
The model for subject $i$ may be:
$$\lambda^{\theta}_{1i}(t;T^i_1,T^i_2,T^i_3)=\alpha_{01}(t)1_{\{T^i_{1}\ge t\}}1_{\{T^i_{3}\ge t\}}e^{\eta_2^1 1_{T^i_2<t}+\gamma_2^1 1_{T^i_2<t}T^i_2+\beta^1 Z_i }$$
$$\lambda^{\theta}_{2i}(t;T^i_1,T^i_2,T^i_3)=\alpha_{02}(t)1_{\{T^i_{2}\ge t\}}1_{\{T^i_{3}\ge t\}}e^{\eta_1^2 1_{T^i_1<t}+\gamma_1^2 1_{T^i_1<t}T^i_1+\beta^2 Z_i };$$
$$\lambda^{\theta}_{3i}(t;T^i_1,T^i_2,T^i_3)=\alpha_{04}(t)1_{\{T^i_{3}\ge t\}}
e^{\eta_1^31_{T^i_1<t}+\eta_2^3 1_{T^i_2<t}+\eta^3_{12}
1_{T^i_1<t}1_{T^i_2<t}+\gamma_2^3 1_{T^i_2<t}T^i_2+\gamma_1^3
1_{T^i_1<t}T^i_1+\beta^3 Z_i }.$$

\subsection{Likelihoods for coarsened observation from the dementia-institution-death model}
Dementia is commonly observed  in discrete time (because it is
diagnosed at planned visits), while time of death and
institutionalization can be known exactly. A heuristic likelihood
for this MDCO scheme of the semi-parametric Markov model was
proposed in Commenges \& Joly (2004). The results of the present
paper allow us to rigorously derive the likelihood, to obtain a more
concise formula and this can also be done for non-Markov models.
The model can be described as a trivariate counting process with
$N_1$ observed at discrete times $v_0, \ldots, v_m$ and $N_2$ and
$N_3$ are observed in continuous time. We will consider that $v_m$ is the last visit really done ($v_m < \tilde T_3$) and thanks to the result of section 4.5 we will be able to treat $v_m$ as deterministic. Thus we have $r_2(t)=r_3(t)=1, 0\le t\le C$ and we observe (for each
subject) $(\tilde T_j, \delta_j),j=2,3$, using the standard notation of section 4.1; note that $\delta_j=1_{C_j}, j=2,3$, with the event $C_j$ (not to be confounded with time $C$) defined in section 4.4.  As for $N_1$ we observe the events $D_{1l}, l=1\ldots,m$ and $E_1$ (or equivalently the value of $1_{D_{1l}}, l=1\ldots,m$ and $1_{E_1}$).
Using the formula of Theorem 2, and
dropping the multiplicative constant, we have for instance on the pseudo-atom $D_{1l}\cap C_2 \cap C_3$:
$$ \LL^{\theta}_{\Ob}= \int_{v_{l-1}}^{v_l}
\lambda_1^{\theta}(s_1;s_1, T_2,
T_3)\lambda_2^{\theta}(T_2;s_1, T_2,
T_3)\lambda_3^{\theta}( T_3;s_1, T_2,
T_3)e^{-\Lambda_.^{\theta}(T_3;s_1,
T_2, T_3)}ds_1. $$
For a more compact expression of the likelihood we may group the formulae for the pseudo-atoms included in $D_{1l}$ by making use of the $\delta_j$'s and the $\tilde T_j$; the likelihood on $D_{1l}=\{T_1\in (v_{l-1},v_l]\}$ is:
\begin{equation} \label{LO1} \LL^{\theta}_{\Ob}= \int_{v_{l-1}}^{v_l}
\lambda_1^{\theta}(s_1;s_1,\tilde T_2,\tilde
T_3)\lambda_2^{\theta}(\tilde T_2;s_1, \tilde T_2,\tilde
T_3)^{\delta_2}\lambda_3^{\theta}(\tilde T_3;s_1,\tilde T_2,\tilde
T_3)^{\delta_3}e^{-\Lambda_.^{\theta}(\tilde T_3;s_1,\tilde
T_2,\tilde T_3)}ds_1. \end{equation}
Similarly we obtain on $E_1=\{T_1>v_m\}$:
\begin{eqnarray} \label{LO2} \LL^{\theta}_{\Ob}&=& \int_{v_{m}}^{\tilde T_3}
\lambda_1^{\theta}(s_1;s_1,\tilde T_2, \tilde
T_3)\lambda_2^{\theta}(\tilde T_2;s_1,\tilde T_2, \tilde
T_3)^{\delta_2}\lambda_3^{\theta}(\tilde T_3;s_1,\tilde T_2, \tilde
T_3)^{\delta_3}e^{-\Lambda_.^{\theta}(\tilde T_3;s_1,\tilde T_2,
\tilde T_3)}ds_1 \nonumber \\ &+&\lambda_2^{\theta}(\tilde
T_2;\tilde T_3,\tilde T_2, \tilde
T_3)^{\delta_2}\lambda_3^{\theta}(\tilde T_3;\tilde T_3,\tilde
T_2, \tilde T_3)^{\delta_3}e^{-\Lambda_.^{\theta}(\tilde
T_3;\tilde T_3,\tilde T_2, \tilde T_3)}, \end{eqnarray} where
$\Lambda_.^{\theta}(t;s,v,
w)=\sum_{j=1}^3\Lambda_j^{\theta}(t;s,v, w)$, and
$\Lambda_j^{\theta}(t;s,v, w)=\int_0^{t}\lambda_j^{\theta}(u;s,v,
w )du$, $j=1,2,3$. In this formula we have replaced the upper bound $C$ of the integral by $\tilde T_3$ because $\lambda_1^{\theta}(s_1;s_1,\tilde T_2, \tilde
T_3)=0, T_3<s_1 \le C$. For the same reason we have replaced $C$ in $\Lambda_.^{\theta}(.;.,.,.,.)$ by $\tilde T_3$.

 For the Markov model (\ref{Markov}) for
instance the intensities appearing inside the integrals of
(\ref{LO1}) and (\ref{LO2}) can be computed to be:
$$\lambda^{\theta}_1(s_1;s_1,\tilde T_2,\tilde T_3)=\alpha_{01}(s_1)e^{\eta_2^1 1_{\{\tilde T_2 < s_1\}}}$$
$$\lambda^{\theta}_2(\tilde T_2;s_1, \tilde T_2,\tilde T_3)=\alpha_{02}(\tilde T_2)e^{\eta_1^2 1_{\{s_1 < \tilde T_2\}}}$$
$$\lambda^{\theta}_3(\tilde T_3;s_1,\tilde T_2,\tilde T_3)=\alpha_{04}(\tilde T_3)e^{\eta_1^3 +\eta_2^3 1_{\{\tilde T_2< \tilde T_3 \}}+\eta_{12}^3 1_{\{\tilde T_2< \tilde T_3 \}}}.$$
Note that the last equality is available only on the set
$\{s_1<\tilde T_3\}$ (which is the case inside the integrals) and note
moreover that on this event $\lambda^{\theta}_3(\tilde
T_3;s_1,\tilde T_2,\tilde T_3)$ does not depend on $s_1$: this comes
from the Markov property of this model. The other quantities
appearing in (\ref{LO1}) and (\ref{LO2}) can be computed by the
same mechanical manipulations.
 The general formula can be applied even if the
model is non-Markovian.

A closer examination of the data available in the PAQUID study has
revealed that the observations were more incomplete that what was assumed
in Commenges \& Joly (2004). When a subject is visited at time
$v_l$ it is observed whether he lives in institution or not and
the time of entry in institution is retrospectively recorded;
information about institutionalization between the last visit and
$\tilde T_3$ is often unknown for subjects who where not
institutionalized at the last visit. The response processes are: $R_{1t}=r_1(t) 1_{\{N_{3t}=0\}}$  where $r_1(t)=1$ for $t=v_l$
, $j=1,\ldots,m_1$ and zero elsewhere. We define  $M$ as
$M=\max_l (l:v_l<\tilde T_3)$, so that $v_M$ is the last visit time before $C$ or death; $R_{2t}=1_{\{t\le v_M\}}$ and $R_{3t}=1$ for $0\le t \le C$. It is
seen that $R_{2t}$ depends on the process $N_3$, and moreover,
because of the retrospective recollection of the time of entry in
institution, it depends on values of $N_3$ for times larger than
$t$ (future values !): so it is not obvious that the mechanism
leading to incomplete data is ignorable. The result of section 4.5 does not apply directly to this case but we can apply  the same approach to prove that ignorability holds in this case.  So we can still
apply the general formula of Theorem 2.

The observed variables for dementia and death are as before, but
for institution we observe $C_2, 1_{C_2}T_2, E_2$. On $C_2\cap D_{1l}$ (that is we observe that the
subject was institutionalized at $T_2$ and demented between
$v_{l-1}$ and $v_l$) the likelihood is given by (\ref{LO1}); on $C_2\cap E_1$ (the
subject was institutionalized at $T_2$ and not demented at the last visit)
 it is given by (\ref{LO2}). Two different formulae are necessary to describe the likelihood on $E_{2}$.
As before the formula of Theorem 2 gives us the likelihoood on pseudo-atoms and we may group the formulae for the pseudo-atoms included in $E_{2}\cap D_{1l}$ to obtain:
\begin{eqnarray} \label{LO3} \LL^{\theta}_{\Ob}&=
&\int_{v_{l-1}}^{v_l} \Bigl [\int_{v_{M}}^{\tilde
T_3}\lambda_1^{\theta}(s_1;s_1,s_2,\tilde
T_3)\lambda_2^{\theta}(s_2;s_1,s_2,\tilde
T_3)\lambda_3^{\theta}(\tilde T_3;s_1,s_2,\tilde
T_3)^{\delta_3}e^{-\Lambda_.^{\theta}(\tilde T_3;s_1,s_2,\tilde
T_3)}ds_2 \nonumber \\ &+&\lambda_1^{\theta}(s_1;s_1,\tilde
T_3,\tilde T_3)\lambda_3^{\theta}(\tilde T_3;s_1,\tilde T_3,\tilde
T_3)^{\delta_3}e^{-\Lambda_.^{\theta}(\tilde T_3;s_1,\tilde
T_3,\tilde T_3)}\Bigr ]ds_1, \end{eqnarray}
 and those incuded in $E_{2}\cap E_1$ to obtain:
\begin{eqnarray} \label{LO4} \LL^{\theta}_{\Ob}&= &\int_{v_{M}}^{\tilde
T_3} \int_{v_{M}}^{\tilde
T_3}\lambda_1^{\theta}(s_1;s_1,s_2,\tilde
T_3)\lambda_2^{\theta}(s_2;s_1,s_2,\tilde
T_3)\lambda_3^{\theta}(\tilde T_3;s_1,s_2,\tilde
T_3)^{\delta_3}e^{-\Lambda_.^{\theta}(\tilde T_3;s_1,s_2,\tilde
T_3)}ds_1ds_2 \nonumber \\ &+&\int_{v_{M}}^{\tilde
T_3}\lambda_1^{\theta}(s_1;s_1,\tilde T_3,\tilde
T_3)\lambda_3^{\theta}(\tilde T_3;s_1,\tilde T_3,\tilde
T_3)^{\delta_3}e^{-\Lambda_.^{\theta}(\tilde T_3;s_1,\tilde
T_3,\tilde T_3)}ds_1\nonumber \\ &+&\int_{v_{M}}^{\tilde
T_3}\lambda_2^{\theta}(s_2;\tilde T_3,s_2,\tilde
T_3)\lambda_3^{\theta}(\tilde T_3;\tilde T_3,s_2,\tilde
T_3)^{\delta_3}e^{-\Lambda_.^{\theta}(\tilde T_3;\tilde
T_3,s_2,\tilde T_3)}ds_2\nonumber \\ &+&\lambda_3^{\theta}(\tilde
T_3;\tilde T_3,\tilde T_3,\tilde
T_3)^{\delta_3}e^{-\Lambda_.^{\theta}(\tilde T_3;\tilde T_3,\tilde
T_3,\tilde T_3)}
. \end{eqnarray}

 Once the likelihood is computed different approaches for inference are possible; in particular it was proposed in Commenges \& Joly (2004) to use penalized likelihood, with smoothing coefficients chosen by cross-validation and this method gave very satisfactory results in a simulation study.

\section{Conclusion}
Many multi-state models can be considered as generated by simple events so that a direct representation in terms of counting processes may be more economical; this is particularly the case when the events are not repeated and can be modeled by $0-1$ counting processes. The multi-state point of view however retains its interest in many applications, in particular for reversible models: for instance  epileptic crises, repeated diarrhoea periods or repeated hospitalization stays might be modeled by a reversible two-state model if we wish to take into account the duration of the crises or of the hospital stays. So the multi-state and the counting process points of view are rather complementary.

 Representing multi-state models as counting processes models allows us to rigorously derive the likelihood by the use of Jacod's formula. This was already known for continuous time observation schemes but we were able to apply this approach to the quite general GCMP scheme. This theory will be useful for developing complex models in life history events analysis. In addition to that, having a compact general formula for the likelihood could be exploited, for instance,  for designing a software able to automatically treat any IM model: the user could specify his model by a routine giving the the values of the intensities of the OJC model which would be used to compute the likelihood corresponding to the observations described in the data set. This would be in particular feasible for parametric or penalized likelihood approaches; see Commenges {\em et al.} (2006) for a penalized likelihood approach of Markov and semi-Markov models.

\subsection*{Acknowledgment:} We thank
Pierre Joly and Alioum Ahmadou for their comments about the
manuscript. We also thank a referee for his careful reading of the paper.

\section{References}
\setlength{\parindent}{0.0in}
\noindent Aalen, O. (1978). Nonparametric inference for a family
of counting processes. {\em Ann. Statist.} {\bf 6},
701-726.\vspace{3mm}

\noindent Aalen, O., Fosen, J., Weedon-Fekjaer, H., Borgan, O. \& Husebye, E. (2004). Dynamic analysis of multivariate failure time data. {\em Biometrics.} {\bf 60},
764-773.\vspace{3mm}

Aalen, O. \& Johansen S. (1978). An empirical
transition matrix for nonhomogenous Markov chains based on
censored observations. {\em  Scand. J. Statist.}
{\bf 5}, 141--150.\vspace{3mm}

Alioum, A.,  Commenges, D., Thi\'ebaut,  R. \& Dabis, F. (2005).  A multi-state approach for estimating incidence of human immunodeficiency virus from prevalent cohort data. {\em Applied Statistics}, {\bf 54}, 739-752.\vspace{3mm}

Andersen, P.K. (1988). Multi-State
models in survival analysis: a study of nephropathy and mortality
in diabetes. {\em Statist. Med.} {\bf 7},
661-670.\vspace{3mm}

Andersen, P.K, Borgan \O, Gill RD
\&  Keiding N. (1993). {\em Statistical Models Based on Counting
Processes}. Springer, New-York.\vspace{3mm}

Borgan, O. (1984). Maximum likelihood estimation in parametric
counting process models, with application to censored failure time
data. {\em Scand. J. Statist.} {\bf 11}, 1-16.\vspace{3mm}

Commenges, D. (2002). Inference for multi-state models
from interval-censored data. {\it Stat. Methods Med.
Res.}, {\bf 11}, 167-182.\vspace{3mm}

Commenges, D. (2003). Likelihood for interval-censored observations from
multi-state models.  {\em Statistics and Operational Research
Transactions} {\bf 27}, 1-12.\vspace{3mm}

 \noindent Commenges, D. \& Joly, P. (2004). Multi-state model for dementia,
institutionalization and death. {\em Commun. Statist.
A} {\bf  33}, 1315-1326.\vspace{3mm}

 \noindent Commenges, D., Joly, P., Letenneur, L. \& Dartigues, JF. (2004). Incidence and
prevalence of Alzheimer's disease or dementia using an
Illness-death model. {\em Statist. Med.} {\bf  23},
199-210.\vspace{3mm}

 \noindent Commenges, D. \& G\'egout-Petit, A. (2005). Likelihood inference for incompletely observed stochastic processes: general ignorability conditions. {\em arXiv:math.ST/0507151}.\vspace{3mm}

 \noindent Commenges, D., Joly, P., G\'egout-Petit, A. \& Liquet, B. (2006). Choice between semi-parametric estimators of Markov and non-Markov multi-state models from generally coarsened observations. {\em Submitted}.\vspace{3mm}

 \noindent Fleming, T. R. (1978).
Nonparametric estimation for nonhomogeneous Markov processes in
the problem of competing risks. {\em Ann. Statist.} {\bf 6},
1057-1070. \vspace{3mm}

 \noindent Frydman, H. (1995a).
Non-parametric estimation of a  Markov ``illness-death model''
process  from interval-censored observations, with application to
diabetes survival data. {\it Biometrika} {\bf  82},
773-789.\vspace{3mm}

 Frydman, H. (1995b). Semi-parametric
estimation in a three-state duration-dependent Markov model from
interval-censored observations with application to AIDS. {\em
Biometrics} {\bf 51}, 502-511. \vspace{3mm}

 Gill, R. D., van der Laan, M. J. \& Robins, J.M. (1997). Coarsening at random:
characterizations, conjectures and counter-examples, in: {\em
State of the Art in Survival Analysis}, D.-Y. Lin \& T.R. Fleming
(eds), Springer Lecture Notes in Statistics 123, 255-294
\vspace{3mm}

 Hougaard, P. (2000). {\em Analysis of
multivariate survival data}, Springer, New York. \vspace{3mm}

\noindent Jacod, J. (1975). Multivariate point processes:
predictable projection; Radon-Nikodym derivative, representation
of martingales. {\em Z. Wahrsch. verw. Geb.} {\bf 31},
235-253.\vspace{3mm}

 \noindent Jacod, J. \& Memin, J. (1976).
Caract\'eristiques locales et conditions de continuit\'e pour les
semi-martingales. {\em  Z. Wahrsch. verw. Geb.} {\bf 35},
1-37.\vspace{3mm}

\noindent Joly, P. \& Commenges, D. (1999). A
penalized likelihood approach for a progressive three-state model
with censored and truncated data: Application to AIDS. {\em
Biometrics} {\bf 55}, 887-890.\vspace{3mm}

 \noindent Joly, P.,
Commenges, D., Helmer, C. \& Letenneur, L. (2002). A penalized
likelihood approach for an illness-death model with
interval-censored data: application to age-specific incidence of
dementia. {\em Biostatistics} {\bf 3}, 433- 443.\vspace{3mm}

Kallenberg, O. (2001). {\em Foundations of modern probabilities}.
Springer, New-York.\vspace{3mm}

 Keiding, N. (1991).
Age-specific incidence and prevalence: a statistical perspective.
{\em  J. Roy. Statist. Soc. Ser. A} {\bf
154}, 371-412.\vspace{3mm}

 Lagakos, S.W., Sommer, C.J. \&
Zelen, M. (1978). Semi-Markov models for partially censored data.
{\em Biometrika} {\bf 65}, 311-317.\vspace{3mm}

 Letenneur, L.,
Gilleron, V., Commenges, D., Helmer, C., Orgogozo, JM. \&
Dartigues JF. (1999). Are sex and educational level independent
predictors of dementia and Alzheimer's disease ? Incidence data
from the PAQUID project. {\em  Journal of Neurology Neurosurgery
and Psychiatry} {\bf 66}, 177-183.\vspace{3mm}

 \noindent  Peto,
R. (1973). Experimental survival curves for interval-censored data
{\em Applied Statistics} {\bf  22}, 86-91.\vspace{3mm}

 \noindent  Schweder, T. (1970). Composable Markov processes. {\em J. Appl. Probab.}
 {\bf 7}, 400-410.\vspace{3mm}

 \noindent Turnbull, B. W. (1976). The
empirical distribution function with arbitrarily grouped, censored
and truncated data. {\it J. Roy. Statist. Soc. Ser. B} {\bf 38}, 290-295.\vspace{3mm}

 \noindent Wong,
G. Y. C. \& Yu Q. (1999). Generalized MLE of a Joint Distribution
Function with Multivariate Interval-Censored Data. {\em J. Multivariate Anal.} {\bf  69}, 155-166.\vspace{5mm}

Daniel Commenges, INSERM E0338; Universit\'e Victor Segalen Bordeaux 2, 146 rue L\'eo
Saignat, Bordeaux, 33076, France\\
E-mail: daniel.commenges@isped.u-bordeaux2.fr

\newpage
 %\vspace{0.5cm}
%\begin{center} %\begin{picture}(350,240)
%\put(-10,100){\framebox(80,40){0: Healthy}}
%\put(150,10){\framebox(80,40){4: Dead}}
%\put(150,100){\framebox(80,40){2: Instit}}
%\put(150,190){\framebox(80,40){1: Ill}}
%\put(310,100){\framebox(80,40){3: Ill+Instit}}
%\put(95,175){\makebox{$\alpha_{01}$}}
%\put(95,125){\makebox{$\alpha_{02}$}}
%\put(95,55){\makebox{$\alpha_{04}$}}
%\put(170,75){\makebox{$\alpha_{23}$}}
%\put(260,175){\makebox{$\alpha_{14}$}}
%\put(260,125){\makebox{$\alpha_{24}$}}
%\put(265,55){\makebox{$\alpha_{34}$}} %\thicklines
%\put(72,120){\vector(1,0){76}} %\put(232,120){\vector(1,0){76}}
%\put(75,95){\vector(1,-1){60}} %\put(75,135){\vector(1,1){60}}
%\put(190,100){\vector(0,-1){48}}
%\put(235,200){\vector(1,-1){60}}
%\put(300,95){\vector(-1,-1){60}} %\put(230,210){\line(1,0){180}}
%\put(410,210){\line(0,-1){180}} %\put(410,30){\vector(-1,0){180}}
%\end{picture}
%\end{center}

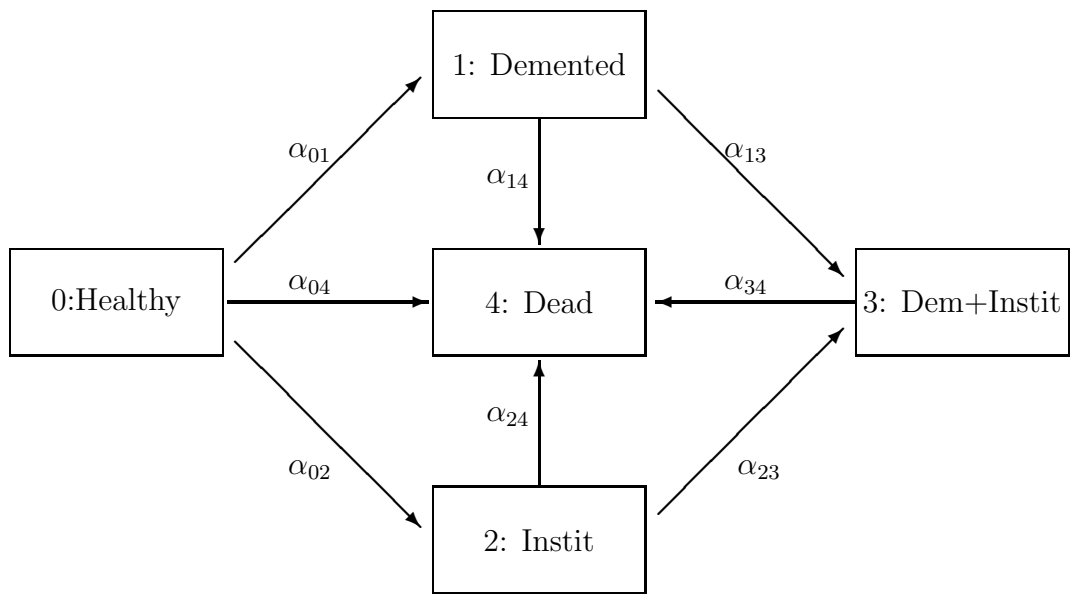
\begin{figure}
\vspace{0.5cm}
\begin{center}
\begin{picture}(350,240)
\put(-10,100){\framebox(80,40){0:Healthy}}
\put(150,10){\framebox(80,40){2: Instit}}
\put(150,100){\framebox(80,40){4: Dead}}
\put(150,190){\framebox(80,40){1: Demented}}
\put(310,100){\framebox(80,40){3: Dem+Instit}}
\put(95,175){\makebox{$\alpha_{01}$}}
\put(95,125){\makebox{$\alpha_{04}$}}
\put(95,55){\makebox{$\alpha_{02}$}}
\put(170,165){\makebox{$\alpha_{14}$}}
\put(170,75){\makebox{$\alpha_{24}$}}
\put(260,175){\makebox{$\alpha_{13}$}}
\put(260,125){\makebox{$\alpha_{34}$}}
\put(265,55){\makebox{$\alpha_{23}$}} \thicklines
\put(72,120){\vector(1,0){76}} \put(310,120){\vector(-1,0){76}}
\put(75,105){\vector(1,-1){70}} \put(75,135){\vector(1,1){70}}
\put(190,190){\vector(0,-1){48}} \put(190,50){\vector(0,1){48}}
\put(235,200){\vector(1,-1){70}} %\put(300,95){\vector(-1,-1){60}}
\put(235,40){\vector(1,1){70}}
\end{picture}
\end{center}
\caption{The five-state model for dementia, institutionalization
and death.}
\end{figure}
\end{document}